\overfullrule=0pt
\centerline {\bf A property shared by continuous linear functions and holomorphic functions}
\bigskip
\bigskip
\centerline {BIAGIO RICCERI}\par
\bigskip
\bigskip
{\bf Abstract:} In this note, we continue to highlight some applications of Theorem 1 of [3]. Here is a
sample: Let $X$ be an open set in ${\bf C}^n$, $\Omega$ an open convex set in ${\bf C}$ and $f, g : X\to
{\bf C}$ two holomorphic functions such that $f(X)\cap\Omega\neq\emptyset$, $f(X)\setminus\Omega\neq
\emptyset$ and $g(X)\subseteq \Omega$. Then, there exists a set $A$ in $[0,1]$ with the following properties:\par\noindent
$(a)$\hskip 5pt for each $x\in X$, there exists $\lambda\in A$ such that $\lambda g(x)+(1-\lambda)f(x)\in\Omega$\ ;\par\noindent 
$(b)$\hskip 5pt for each finite set $B$ in $A$, there exists $u\in X$ such that
$\mu g(u)+(1-\mu)f(u)\in {\bf C}\setminus\Omega$ for all $\mu\in B$.\par
 \par
\bigskip
\bigskip
\bigskip
\bigskip
This paper is essentially a companion to [3].\par
\smallskip
In the sequel, the term "interval" means a non-empty connected subset of ${\bf R}$ with more than one point.\par
\smallskip
  For 
a multifunction $F:I\to 2^X$, as usual, for $A\subseteq I$ and $B\subseteq X$, we set
$$F(A)=\bigcup_{x\in A}F(x)\ ,$$
$$F^{-}(B)=\{\lambda\in I : F(\lambda)\cap B\neq\emptyset\}$$
and
$$F^{+}(B)=\{\lambda\in I : F(\lambda)\subseteq B\}\ .$$
When $I$ is an interval,  $F$ 
  is said to be non-decreasing (resp.
non-increasing) with respect to the inclusion if $F(\lambda)\subseteq F(\mu)$ 
(resp. $F(\mu)\subseteq F(\lambda))$ for all $\lambda, \mu\in I$, with $\lambda<\mu$.\par
\medskip
Furthermore, let $Y$ be a non-empty set and ${\cal F}$ a family of
subsets of $Y$. We say that ${\cal F}$ has the compactness-like property if
every subfamily of ${\cal F}$ satisfying the finite intersection property has
a non-empty intersection.  An obvious (but useful) remark is that if ${\cal F}$
has the compactness-like property and $\Gamma_0\neq\emptyset$ is a member of
${\cal F}$, then the family $\{\Gamma\cap \Gamma_0\}_{\Gamma\in {\cal F}}$ has
the compactness-like property too.
\smallskip
In [3], we established the following general result:\par
\medskip
THEOREM A. - {\it Let $X$ be a non-empty set, $I\subseteq {\bf R}$ an
 interval and $F:I\to 2^X$ a multifunction such that,
for each $x\in X$, the set $X\setminus F^-(x)$ is an interval
open in $I$. Moreover, assume that, for some $\lambda_0\in I$, with
$F(\lambda_0)\neq\emptyset$, and for some set $D\subseteq I$ dense in $I$, the family $\{F(\lambda)\cap F(\lambda_0)\}_{\lambda\in D}$
has the compactness-like property.\par
Under such hypotheses, there exists a compact interval
$[a^*,b^*]\subseteq I$  such that either $a^*\geq \lambda_0$, 
$( F(a^*)\cap F(\lambda_0))\setminus F(]a^*,b^*[)\neq \emptyset$ and $F_{|]a^*,b^*[}$
is non-decreasing with respect to the inclusion, or $b^*\leq \lambda_0$,
$( F(b^*)\cap F(\lambda_0))\setminus F(]a^*,b^*[)\neq \emptyset$ and $F_{|]a^*,b^*[}$  
is non-increasing with respect to the inclusion. 
In particular, the first (resp. second)
occurrence is true when $\lambda_0=\inf I$ (resp. $\lambda_0=\sup I$). Furthermore, if, for some 
neighbourhood
$U$ of $\lambda_0$ in $I$, one has $\cap_{\lambda\in U}F(\lambda)\neq\emptyset$, then in the first
(resp. second) occurrence one has $a^*>\lambda_0$ (resp. $b^*<\lambda_0$).}\par
\medskip
REMARK 1. - We want to remark that the very final part of the conclusion is not present in the
formulation given in [1]. Such a further information follows immediately from the proof. Hence, we refer
the reader to [3] and limit ourselves here to highlight the relevant point only, keeping the notations of the proof
in [1]. So, assume that, for some $\delta>0$,
 one has $\cap_{\lambda\in [\lambda_0-\delta,\lambda_0+\delta]\cap I}F(\lambda)\neq\emptyset$. 
Pick $\tilde x\in \cap_{\lambda\in [\lambda_0-\delta,\lambda_0+\delta]\cap I}F(\lambda)$. Hence,
either $\alpha(\tilde x)\geq \lambda_0+\delta$ or $\beta(\tilde x)\leq \lambda_0-\delta$. If 
$\alpha(\tilde x)\geq \lambda_0+\delta$ (resp. $\beta(\tilde x)\leq \lambda_0-\delta$),
we know that
the first (resp. second) occurrence of the conclusion is true with $a^*=\sup_X\alpha$ (resp. $b^*=\inf_X\beta$),
and the claim follows.\par
\medskip
Within the project of a systematic use of Theorem A, we obtained, in particular, Theorem 3 of [3]. \par
\smallskip
Our aim in the present paper is to give an extension of this latter result to multifunctions (Theorem 2) and then
an application which highlights the property mentioned in the title (Theorem 3).
\smallskip
In the sequel, $Y$ is a real or complex Hausdorff locally convex topological vector space, $E$ is a
closed subset of $Y$ such that $Y\setminus E$ is convex, $I\subseteq {\bf R}$ is an interval,
$\lambda_0\in I$ and $G:X\times I\to 2^Y$ is a given multifunction. The symbol $\partial$ stands
for boundary.\par
\smallskip
We recall a further notion on multifunctions and refer to [1], [2] for other classical notions.\par
\smallskip
Namely, a multifunction $\Phi:I\to 2^Y$ is said to be concave if one has
$$\Phi(t\lambda+(1-t)\mu)\subseteq t\Phi(\lambda)+(1-t)\Phi(\mu)$$
for all $t\in [0,1]$, $\lambda,\mu\in I$.\par
\smallskip
Our main result is as follows:\par
\medskip
THEOREM 1. - {\it  Let the following assumptions hold:\par
\noindent
$(i_1)$\hskip 5pt for each $x\in X$, the multifunction  
$G(x,\cdot)$
is concave and continuous in $I$, and $G(x,\lambda)\subseteq Y\setminus E$
for some $\lambda\in I$\ ;\par
\noindent
$(i_2)$\hskip 5pt 
the set $\{x\in X : G(x,\lambda_0)\cap \hbox {\rm int}(E)\neq\emptyset\}$ is non-empty\ ;\par
\noindent
$(i_3)$\hskip 5pt there exists a set $D\subseteq I$ dense in $I$ such that,
  for every set $A\subseteq D$ for which 
$$\bigcap_{\lambda\in B\cup \{\lambda_0\}}\{x\in X : G(x,\lambda)\cap E\neq\emptyset\}\neq\emptyset$$ for each finite set
$B\subseteq A$, one has
$$\bigcap_{\lambda\in A\cup \{\lambda_0\}}\{x\in X : G(x,\lambda)\cap E\neq\emptyset\}\neq\emptyset\ .$$\par
Then, there exist a compact interval $[a^*,b^*]\subseteq I$ and a point $x^*\in X$, with
$G(x^*,\lambda_0)\cap E\neq\emptyset$ and 
$G(x^*,]a^*,b^*[)\subseteq Y\setminus E$, such that, if we put
$$V=\bigcup_{\lambda\in ]a^*,b^*[}\{x\in X : G(x,\lambda)\subseteq Y\setminus E\}\ ,$$
at least one of the following holds:\par
\noindent
$(p_1)$\hskip 5pt $a^*>\lambda_0$, $G(x^*,a^*)\cap E\neq\emptyset$ and
$$G(V,a^*)\cap E\subseteq \partial G(V,a^*)\cap\partial E\ ;$$
$(p_2)$\hskip 5pt $b^*<\lambda_0$, $G(x^*,b^*)\cap E\neq\emptyset$ and
$$G(V,b^*)\cap E\subseteq \partial G(V,b^*)\cap\partial E\ .$$
In particular, $(p_1)$ (resp. $(p_2)$) holds when $\lambda_0=\inf I$ (resp. $\lambda_0=\sup I$).}\par
\smallskip
PROOF.  Consider the multifunction $F:I\to 2^X$ defined by
$$F(\lambda)=\{x\in X : G(x,\lambda)\cap E\neq\emptyset\}$$
for all $\lambda\in I$.  Note that, by $(i_2)$, there is $\tilde x\in X$ such that $G(\tilde x,\lambda_0)\cap
\hbox {\rm int}(E)\neq \emptyset$. Then, since $G(\tilde x,\cdot)$ is lower semicontinuous, there is a
neighbourhood $U$ of $\lambda_0$ in $I$ such that $G(\tilde x,\lambda)\cap \hbox {\rm int}(E)\neq \emptyset$
for all $\lambda\in U$. Therefore, $\tilde x\in \cap_{\lambda\in U}F(\lambda)$.
In view of $(i_3)$,  the family $\{F(\lambda)\cap
F(\lambda_0)\}_{\lambda\in D}$ has the compactness-like property. Moreover, by $(i_1)$, for each $x\in X$, the set
$I\setminus F^{-}(x)$ (that is $\{\lambda\in I : G(x,\lambda)\subseteq Y\setminus E\}$) is non-empty,
convex (since $G(x,\cdot)$ is concave and $Y\setminus E$ is convex) and open in $I$ (since $G(x,\cdot)$
is upper semicontinuous and $Y\setminus E$ is open). Hence, the multifunction $F$ satisfies the hypotheses of
Theorem A. Consequently, there exists a compact interval $[a^*,b^*]\subseteq I$ such that either
$a^*>\lambda_0$, 
$( F(a^*)\cap F(\lambda_0))\setminus F(]a^*,b^*[)\neq \emptyset$ and $F_{|]a^*,b^*[}$
is non-decreasing with respect to the inclusion, or $b^*<\lambda_0$,
$( F(b^*)\cap F(\lambda_0))\setminus F(]a^*,b^*[)\neq \emptyset$ and $F_{|]a^*,b^*[}$  
is non-increasing with respect to the inclusion. For instance, assume that the first alternative holds.
Pick $x^*\in ( F(a^*)\cap F(\lambda_0))\setminus F(]a^*,b^*[)$. So, 
$G(x^*,\lambda_0)\cap E\neq\emptyset$, $G(x^*,a^*)\cap E\neq\emptyset$ and 
$G(x^*,]a^*,b^*[)\subseteq Y\setminus E$. First,  let us prove that $G(V,a^*)\cap E\subseteq \partial E$.
So, let  $y\in G(V,a^*)\cap E$.
Arguing by contradiction, assume that $y\in \hbox {\rm int}(E)$. Let $x\in
V$ be such that $y\in G(x,a^*)$. Also, let $\mu\in ]a^*,b^*[$ be such that
$G(x,\mu)\subseteq Y\setminus E$. Since $G(x,a^*)\cap \hbox {\rm int}(E)\neq\emptyset$
and $G(x,\cdot)$ is lower semicontinuous, we can find $\lambda\in ]a^*,\mu[$ so that
$G(x,\lambda)\cap \hbox {\rm int}(E)\neq\emptyset$. But then $X\setminus F(\mu)\subseteq 
X\setminus F(\lambda)$ and so $G(x,\lambda)\subseteq Y\setminus E$, a contradiction. Next,
let us prove that $G(V,a^*)\cap E\subseteq \partial G(V,a^*)$. Fix $z\in G(V,a^*)\cap E$.
Arguing by contradiction again, assume that $z\in \hbox {\rm int}(G(V,a^*))$. We already know
that $z\in \partial E$ and so $z\in \partial(Y\setminus E)$. Since $Y\setminus E$ is open and convex,
a classical separation theorem ensures the existence of  a non-zero continuous linear
functional $\varphi : Y\to {\bf R}$ such that
$\varphi(z)<\varphi(u)$ for all $u\in Y\setminus E$. Therefore, the set $\varphi^{-1}(]-\infty,\varphi(z)[)$
is contained in $E$. Moreover, it is open and meets $\hbox {\rm int}(G(V,a^*))$ since, otherwise,
$z$ would be a local minimum of $\varphi$, which is impossible since $\varphi$ is linear. Consequently, the 
set $\varphi(]-\infty,\varphi(z)[)\cap \hbox {\rm int}(G(V,a^*))$ would be non-empty, open and containted
$G(V,a^*)\cap E$. So, by what we have already seen, we would have $\varphi(]-\infty,\varphi(z)[)\cap \hbox {\rm int}(G(V,a^*))\subseteq
\partial E$, against the fact that int$(\partial E)=\emptyset$ since $E$ is closed. So $(p_1)$ is proved.
If the second alternative above holds, we obtain $(p_2)$ by means of analogous reasonings.\hfill
$\bigtriangleup$
\medskip
REMARK 2. - Notice that the mere lower semicontinuity of $G(x,\cdot)$ is not enough for the validity of Theorem 1. To see this,
take $Y={\bf R}$, $E={\bf R}\setminus ]-1,1[$, $I=[0,1]$ and $G(x,\lambda)=]-1,1[-\lambda$ for all $(x,\lambda)\in X\times
[0,1]$. Clearly, each assumptions of Theorem 1, but the upper semicontinuity of $G(x,\cdot)$, is satisfied. However, neither
$(p_1)$ nor $(p_2)$ holds since the values of $G$ are open.\par
\medskip
REMARK 3. - Clearly, if we give up the better information $a^*>\lambda_0$ (resp. $b^*<\lambda_0$), the continuity of
$G(x,\cdot)$ can be weakened to upper semicontinuity and, at the same time,  $(i_2)$ can be weakened to
 $\{x\in X : G(x,\lambda_0)\cap E\neq\emptyset\}\neq\emptyset$.\par
\medskip
The above-mentioned extension of Theorem 3 of [3] obtained via Theorem 1 is as follows:\par
\medskip
THEOREM 2. - {\it Let $0\in I$ and let $\Phi, \Psi$ be two multifunctions from $X$ into $Y$, with non-empty compact values. Assume that:\par
\noindent
$(h_1)$\hskip 5pt for each $x\in $X, there exists $\lambda\in I$ such that $\Phi(x)+\lambda\Psi(x)\subseteq Y\setminus E$\ ;\par
\noindent
$(h_2)$\hskip 5pt the set $\Phi^{-}(\hbox {\rm int}(E))$ is non-empty\ ;\par
\noindent
$(h_3)$\hskip 5pt there exists a set $D\subseteq I$ dense in $I$ such that, for every set
$A\subseteq D$ for which 
$$\bigcap_{\lambda\in  B\cup \{0\}}(\Phi+\lambda\Psi)^{-}(E)\neq\emptyset$$ for each finite set
$B\subseteq A$, one has
$$\bigcap_{\lambda\in A\cup \{0\}}(\Phi+\lambda\Psi)^{-}(E)\neq\emptyset\ .$$
Then, there exist a compact interval $[a^*,b^*]\subseteq I$ and a point $x^*\in \Phi^{-}(E)$, with
$$\bigcup_{\lambda\in ]a^*,b^*[}(\Phi(x^*)+\lambda\Psi(x^*))\subseteq Y\setminus E\ ,$$ such that, if we put
$$V=\bigcup_{\lambda\in ]a^*,b^*[}(\Phi+\lambda\Psi)^{+}(Y\setminus E)\ ,$$
at least one of the following holds:\par
\noindent
$(q_1)$\hskip 5pt $a^*>0$, $(\Phi(x^*)+a^*\Psi(x^*))\cap E\neq\emptyset$ and
$$(\Phi+a^*\Psi)(V)\cap E\subseteq \partial (\Phi+a^*\Psi)(V)\cap\partial E\ ;$$
$(q_2)$\hskip 5pt $b^*<0$, $(\Phi+b^*\Psi)(V)\cap E\neq\emptyset$ and
$$(\Phi+b^*\Psi)(V)\cap E\subseteq \partial (\Phi+b^*\Psi)(V)\cap\partial E\ .$$
In particular, $(q_1)$ (resp. $(q_2)$)
holds when $0=\inf I$ (resp. $0=\sup I$)\ .}\par
\smallskip
PROOF. Apply Theorem 1, with $\lambda_0=0$, to the multifunction
 $G:X\times I\to 2^Y$ defined by
$$G(x,\lambda)=\Phi(x)+\lambda\Psi(x)$$
for all $(x,\lambda)\in X\times I$. In particular, note that the multifunction $G(x,\cdot)$ is continuous
 since both $\Phi(x), \Psi(x)$ are compact. \hfill $\bigtriangleup$\par
\medskip
REMARK 4. - According to Remark 3, if, instead of $(h_2)$, we simply assume that 
$\Phi^{-}(E)\neq\emptyset$, then $(q_1)$ (resp. $(q_2)$) holds with $a^*\geq 0$ (resp.
$b^*\leq 0$).\par
\medskip
Here is the property mentioned in the title:\par
\medskip
THEOREM 3. - {\it Let $f, g :X\to Y$ be two functions such that 
$f(X)\cap E\neq\emptyset$, $f(X)\setminus E\neq\emptyset$ and $g(X)\subseteq Y\setminus E$. Moreover,
suppose that one of the two following sets of assumptions holds:\par
\noindent
$(k_1)$\hskip 5pt $X$ is an open set in a Banach space $S$, $\hbox {\rm dim}(Y)<\infty$
 and
$f, g$ are the restrictions to $X$ of two continuous linear functions from $S$ into $Y$\ ;\par
\noindent
$(k_2)$\hskip 5pt $X$ is an open set in ${\bf C}^n$, $Y={\bf C}$ and $f, g$ are holomorphic in $X$.\par
Then, for every set $D\subseteq [0,1]$ dense in $[0,1]$, with $0\in D$, there exists a set $A\subseteq D$
with the following properties:\par
\noindent
$(r_1)$\hskip 5pt for every $x\in X$, there exists $\lambda\in A$ such that
$$\lambda g(x)+(1-\lambda)f(x)\in Y\setminus E\ ;$$
$(r_2)$\hskip 5pt for every finite set $B\subseteq A$, there exists $u\in X$ such that
$$\mu g(u)+(1-\mu)f(u)\in E$$
for all $\mu\in B$.}\par
\smallskip
PROOF. Arguing by contradiction, assume that the conclusion is false. So, assume that there is
a set $D\subseteq [0,1]$ dense in $[0,1]$, with $0\in D$, such that, for every set
$A\subseteq D$ for which 
$$\bigcap_{\lambda\in  B}(\lambda g+(1-\lambda) f)^{-1}(E)\neq\emptyset$$ for each finite set
$B\subseteq A$, one has
$$\bigcap_{\lambda\in A}(\lambda g+(1-\lambda) f)^{-1}(E)\neq\emptyset\ .$$
This means that the family $\{(\lambda g+(1-\lambda)f)^{-1}(E)\}_{\lambda\in D}$ has the compactness-like
property. Then, since $0\in D$, the family $\{(\lambda g+(1-\lambda)f)^{-1}(E)\cap f^{-1}(E)\}_{\lambda\in D}$
has the compactness-like property too. Hence, if we take $I=[0,1]$, $\Phi=f$ and $\Psi=g-f$, assumption $(h_3)$
of Theorem 2 is satisfied.  Moreover, by assumption, $f^{-1}(E)\neq\emptyset$ and, for each $x\in X$, there is
$\lambda\in [0,1]$ such that $f(x)+\lambda(g(x)-f(x))\in Y\setminus E$) (actually, we can take $\lambda=0$ if
$f(x)\in Y\setminus E$, and $\lambda=1$ if $f(x)\in E$). 
Hence, by Theorem 2 (recalling Remark 4), there exist
an interval $[a^*,b^*]\subseteq [0,1]$ and a point $x^*\in f^{-1}(E)$, with
$$w:=a^*g(x^*)+(1-a^*)f(x^*)\in E$$
and
$$\lambda g(x^*)+(1-\lambda) f(x^*)\in Y\setminus E$$
for all $\lambda\in ]a^*,b^*[$,
such that, if we put
$$V=\bigcup_{\lambda\in ]a^*,b^*[}(\lambda g+(1-\lambda) f)^{-1}(Y\setminus E)\ ,$$
we have
$$(a^* g+(1-a^*)f)(V)\cap E\subseteq \partial(a^* g+(1-a^*)f)(V)\cap\partial E\ .\eqno{(1)}$$
Now, assume that $(k_1)$ holds. Continue to denote by the same symbols the continuous linear
extensions of $f, g$ to $S$. Set
$$T=(a^*g+(1-a^*)f)(S)\ .$$
Since $a^*g+(1-a^*)f$ is linear and continuous and dim$(T)<\infty$, by the open mapping theorem,
the set $(a^*g+(1-a^*)f)(V)$ is open in $T$ since $V$ is open in $S$. Then, since $x^*\in V$, there
exists an open neighbourhood $W$ of $w$ in $Y$ such that
$$W\cap T\subseteq (a^*g+(1-a^*)f)(V)\ .$$
From this, it follows that
$$W\cap T\cap E\subseteq (a^*g+(1-a^*)f)(V)\cap E$$
and so, in view of $(1)$,
$$W\cap T\cap E\subseteq \partial E\ .\eqno{(2)}$$ 
Then, since $w\in W\cap T\cap E$, we have $w\in\partial(Y\setminus E)$. Let $\varphi:Y\to {\bf R}$
be a non-zero continuous linear functional such that $\varphi(u)<\varphi(w)$ for
all $u\in Y\setminus E$. By assumption, there is $\tilde x\in X$ such that $f(\tilde x)\in Y\setminus E$.
Set
$$\tilde w=a^*g(\tilde x)+(1-a^*)f(\tilde x)\ .$$
Hence, $\varphi(\tilde w)<\varphi(w)$ as $\tilde w\in Y\setminus E$. Now, fix $\lambda<0$ so that
$w+\lambda(\tilde w-w)\in W$. Then, observing that 
$Y\setminus \overline {Y\setminus E}=$int$(E)$ and that $\varphi(w+\lambda(\tilde w-w))>\varphi(w)$,
we have 
$$w+\lambda(\tilde w-w)\in \hbox {\rm int}(E)$$
which contradicts $(2)$.\par
Now, suppose that $(k_2)$ holds. 
Since $a^*g+(1-a^*)f$ is holomorphic and not constant (note that
$w\in\partial E$ and $\tilde w\in {\bf C}\setminus E$),  and $V$ is open in ${\bf C}^n$,  by a classical
result, the set $(a^*g+(1-a^*)f)(V)$ is open in ${\bf C}$
and this is against $(1)$. The proof is complete.\hfill $\bigtriangleup$\par

\vfill\eject
\centerline {\bf References}\par
\bigskip
\bigskip
\noindent
[1]\hskip 5pt S. HU and N. S. PAPAGEORGIOU, {\it Handbook of multivalued analysis}, Vol. I,  Kluwer Academic Publishers,
1997.\par
\smallskip
\noindent
[2]\hskip 5pt E. KLEIN and A. C. THOMPSON,
{\it Theory of correspondences}, John Wiley \& Sons, 1984. \par
\smallskip
\noindent
[3]\hskip 5pt B. RICCERI, {\it Revisiting a theorem on multifunctions of one real variable},
J. Nonlinear Convex Anal., {\bf 15} (2014), 539-546.\par
\bigskip
\bigskip
\bigskip
\bigskip
Department of Mathematics\par
University of Catania\par
Viale A. Doria 6\par
95125 Catania\par
Italy\par
\smallskip
{\it e-mail address}: ricceri@dmi.unict.it

\bye